\documentclass[12pt]{amsart}
\usepackage{amssymb}

\begin{document}
\bibliographystyle{alpha}
\numberwithin{equation}{section}

\def\Label#1{\label{#1}}

\def\1#1{\ov{#1}}
\def\2#1{\widetilde{#1}}
\def\3#1{\mathcal{#1}}
\def\4#1{\widehat{#1}}

\def\s{s}
\def\k{\kappa}
\def\ov{\overline}
\def\span{\text{\rm span}}
\def\tr{\text{\rm tr}}
\def\GL{{\sf GL}}
\def\xo {{x_0}}
\def\Rk{\text{\rm Rk\,}}
\def\sg{\sigma}

\def \hn {holomorphically nondegenerate}
\def\hyp{hypersurface}
\def\prt#1{{\partial \over\partial #1}}
\def\det{{\text{\rm det}}}
\def\wob{{w\over B(z)}}
\def\co{\chi_1}
\def\po{p_0}
\def\fb {\bar f}
\def\gb {\bar g}
\def\Fb {\ov F}
\def\Gb {\ov G}
\def\Hb {\ov H}
\def\zb {\bar z}
\def\wb {\bar w}
\def \qb {\bar Q}
\def \t {\tau}
\def\z{\chi}
\def\w{\tau}
\def\Z{\zeta}

\def \T {\theta}
\def \Th {\Theta}
\def \L {\Lambda}
\def\b{\beta}
\def\a{\alpha}
\def\o{\omega}
\def\l{\lambda}

\def \im{\text{\rm Im }}
\def \re{\text{\rm Re }}
\def \Char{\text{\rm Char }}
\def \supp{\text{\rm supp }}
\def \codim{\text{\rm codim }}
\def \Ht{\text{\rm ht }}
\def \Dt{\text{\rm dt }}
\def \hO{\widehat{\mathcal O}}
\def \cl{\text{\rm cl }}
\def \bR{\mathbb R}
\def \bC{\mathbb C}
\def \bP{\mathbb P}
\def \C{\mathbb C}
\def \bL{\mathbb L}
\def \bZ{\mathbb Z}
\def \bN{\mathbb N}
\def \scrF{\mathcal F}
\def \scrK{\mathcal K}
\def \scrM{\mathcal M}
\def \cR{\mathcal R}
\def \scrJ{\mathcal J}
\def \scrA{\mathcal A}
\def \scrO{\mathcal O}
\def \scrV{\mathcal V}
\def \scrL{\mathcal L}
\def \scrE{\mathcal E}
\def \hol{\text{\rm hol}}
\def \aut{\text{\rm aut}}
\def \Aut{\text{\rm Aut}}
\def \J{\text{\rm Jac}}
\def\jet#1#2{J^{#1}_{#2}}
\def\gp#1{G^{#1}}
\def\gpo{\gp {2k_0}_0}
\def\emmp {\scrF(M,p;M',p')}
\def\rk{\text{\rm rk}}
\def\Orb{\text{\rm Orb\,}}
\def\Exp{\text{\rm Exp\,}}
\def\ess{\text{\rm Ess\,}}
\def\mult{\text{\rm mult\,}}
\def\Jac{\text{\rm Jac\,}}
\def\Span{\text{\rm span\,}}
\def\d{\partial}
\def\D{\3J}
\def\pr{{\rm pr}}
\def\dbl{[\![}
\def\dbr{]\!]}
\def\nl{|\!|}
\def\nr{|\!|}

\def \depth{\text{\rm depth\,}}
\def \D{\text{\rm Der}\,}
\def \Rk{\text{\rm Rk}\,}
\def \ima{\text{\rm im}\,}
\def \vfi{\varphi}

\title [Transversality of holomorphic mappings between real
hypersurfaces]{Transversality of holomorphic mappings between
real hypersurfaces in  different dimensions}
\author[M.S. Baouendi, P. Ebenfelt, and L. P. Rothschild]{M.S. Baouendi, Peter
Ebenfelt, and Linda P. Rothschild} \thanks{{\rm The first and third
authors are supported in part by the NSF grant DMS-0400880. The
second author is supported in part by the NSF grant DMS-0401215.
\newline}}

\address{ Department of Mathematics, University of California at San Diego, La
Jolla, CA 92093-0112, USA}
\email{sbaouendi@ucsd.edu,  pebenfel@math.ucsd.edu, lrothschild@ucsd.edu }


\abstract In this paper, we consider holomorphic mappings between real hypersurfaces in
different dimensional complex spaces. We give a number of conditions that imply that
such mappings are  transversal to the target hypersurface at most
points.
\endabstract

\newtheorem{Thm}{Theorem}[section]
\newtheorem{Def}[Thm]{Definition}
\newtheorem{Cor}[Thm]{Corollary}
\newtheorem{Pro}[Thm]{Proposition}
\newtheorem{Lem}[Thm]{Lemma}
\newtheorem{Rem}[Thm]{Remark}
\newtheorem{Ex}[Thm]{Example}

\maketitle
\section{Introduction and Main Results}  \Label{s:intro}

In this paper, we consider holomorphic mappings between real hypersurfaces in
different dimensional complex spaces. We give a number of conditions implying
that such  mappings are  transversal to the target hypersurface at
most points. Recall that if
$U$ is an open subset of
$\bC^{n+1}$, $H$ a holomorphic mapping $ U\to \bC^{n'+1}$, and
$M'$ a real hypersurface through a point $H(p)$ for some $p\in U$, then $H$ is
said to be {\it  transversal} to $M'$ at $H(p)$ if
\begin{equation}
T_{H(p)}
M'+dH(T_p\bC^{n+1})=T_{H(p)}\bC^{n'+1},
\end{equation}
where $T_p\C^{n+1}$ and $T_{H(p)}M'$ denote the real tangent
spaces of $\C^{n+1}$ and $M'$ at
$p$ and $H(p)$, respectively.  (We mention that the notion of transversality of a mapping to a hypersurface coincides with that of CR transversality; cf.\ \cite{er06}.) We shall assume that there is a
real hypersurface
$M\subset U$ such that $H(M)\subset M'$. Then
transversality at a point
$H(p)$, for $p\in M$, is equivalent to the nonvanishing at $p$ of the normal derivative
of the real function $u:=\rho'\circ H$, where $\rho'=0$ is a local defining equation for
$M'$ near $H(p)$. Hence a result on  transversality can be
regarded as a type  of Hopf Lemma. 

The equidimensional case (i.e.\
$n=n'$) has been considered by many authors; we mention here the
papers \cite{fornaess76}, \cite{fornaess78}, \cite{p77},
\cite{bb82}, \cite{brgeom}, \cite{brhopf},  \cite{cr94}, \cite{cr98}, \cite{BHR95},
\cite{hp96}, \cite{er06}, \cite{lm06}. In the equidimensional case, 
transversality holds at a point
$H(p)$ under rather general conditions. For instance, in \cite{er06} it is proved that
$H$ is  transversal to $M'$ at $H(p)$ provided that $M'$ is of finite
type at $H(p)$ and the generic rank of $H|_{\Sigma_p}$, where
$\Sigma_p$ denotes the Segre variety of $M$ at $p$, is $n$. The
situation in the case where $n'>n$ is much more complicated. Indeed,
 transversality may fail at a point $H(p)$ even for a polynomial
embedding $\bC^2\to \bC^3$ sending one nondegenerate
hyperquadric into another, as is illustrated by Example \ref{ex2}
below. Observe that a trivial case where  transversality fails at all
points is when $H(U)$ is contained in $M'$. In the equidimensional
case, this is the only way for  transversality to fail at all points
provided that $M$ is holomorphically nondegenerate (see Example
\ref{ex0} and Theorem
\ref{thm5}). When
$n'>n$, this is no longer the case, as is illustrated by Theorem \ref{Thm4} as well
as Example
\ref{ex3}. Our Theorems \ref{main1} and \ref{thm3} give conditions that guarantee
 transversality at most points. The results are essentially optimal, as is
illustrated by examples. Having  tranversality at most points is
crucial  in the study of rigidity of  embeddings into hyperquadrics.
See e.g.\
\cite{w79}, \cite{f86}, \cite{cs83}, \cite{fo86}, \cite{da88}, \cite{huangjdg},
\cite{ehz1}, \cite{ehz2}, \cite{bh}.  See also  \cite{z07} and \cite{a07} for recent related work on transversality of holomorphic Segre mappings.

Before stating our main results, we introduce some notation. 
Let $M$ be a hypersurface in $\bC^{n+1}$, $p\in M$, and $\mathcal L\colon
\bC^{n}\times\bC^{n}\to \bC$ a representative of the Levi form of $M$ at $p$. We
shall denote by $e(M,p):=\min (e_-,e_+)$ and $e_0(M,p)=e_0$, where $e_+,e_-,e_0$,
denote the number of positive, negative, and zero eigenvalues of $\mathcal L$ at $p$.
Observe that $e(M,p)$ and $e_0(M,p)$ are independent of the choice of representative
$\mathcal L$ of the Levi form. A connected hypersurface $M$ is said to be {\it
holomorphically nondegenerate} if there are no germs of nontrivial holomorphic
$(1,0)$-vector fields tangent to $M$. We point out that if $M$ is connected and Levi
nondegenerate at some point, i.e.\ $e_0(M,p)=0$ for some $p\in M$, then $M$ is
necessarily holomorphically nondegenerate. The converse is not true. The reader is
referred to \cite{book} for further details on this and other related notions (see also
\cite{stanton} for holomorphic nondegeneracy). 

In our first theorem, we give two independent conditions
guaranteeing  transversality at most points. 

\begin{Thm} \Label{main1}Let $M\subset \C^{n+1}$, and $M'\subset \C^{n'+1}$ be
connected real-analytic hypersurfaces and $U$ an open neighborhood of $M$
in
$\C^{n+1}$.  Assume that $M$ is holomorphically nondegenerate and that either
\begin{equation}\Label{eigen10} e(M',p') +e_0(M',p') \le n-1,\quad \forall p'\in M'
\end{equation} 
or 
\begin{equation}\Label{eigen20} n' +e_0(M',p') \leq 2n, \quad  \
\ \forall p'\in M',
\end{equation}
 holds.  If $H:U\to \C^{n'+1}$ is a holomorphic
mapping with $H(M) \subset M'$, then one of the following two mutually exclusive 
conditions holds.
\begin{enumerate}
\item [(i)] There is an open subset $V\subset U$ with $M\subset V$ and 
$H(V) \subset
M'$.
\item [(ii)] $H$ is  transversal to $M'$ at $H(p)$ for all $p\in M$
outside a proper real-analytic subset.
\end{enumerate}
\end{Thm}

\begin{Rem} \Label{rmk1}  {\rm We point out that (i) holds if and only if there
exists a point $p\in M$ and a open neighborhood $W\subset U$ of $p$ in
$\bC^{n+1}$ such that $H(W)\subset M'$. This follows easily from the
connectedness and real-analyticity of $M$. Similarly, (ii) holds if and only if
there exists $p\in M$ such that $H$ is  transversal to $M'$ at
$H(p)$.  Indeed, if $H$ is  transversal at $p\in M$, then we have
$\rho'\circ H=a\rho$, where $\rho$ and $\rho'$ are local real-analytic
defining functions near $p$ and
$H(p)$ respectively and $a$ is a real-analytic function defined near
$p$, with $a(p)\neq 0$. The set of points, near $p$, at
which $H$ is not  transversal is given by the equation $a=0$, which
defines a proper real-analytic subset of $M$ near $p$. A standard
connectedness argument shows that (ii) holds. 
 }
\end{Rem}

The condition \eqref{eigen10} in Theorem \ref{main1} is optimal, as can be seen by
Example \ref{ex3}. Similarly, Theorem \ref{Thm4} below shows that condition
\eqref{eigen20} is also optimal. However, if $M$ and $M'$ are Levi nondegenerate and
the target is a hyperquadric\footnote{By a hyperquadric in $\bC^{n+1}$, we
mean a real-algebraic hypersurface defined by $\im w=\langle z,\bar z\rangle$,
where $\langle\cdot,\cdot \rangle$ is a hermitian form in $\bC^n$.}, then condition
\eqref{eigen20} in Theorem
\ref{main1} can be weakened, as is shown by the following result. 

\begin{Thm}\Label{thm3} Let $M\subset \C^{n+1}$ be a connected,
real-analytic hypersurface and $U$ an open neighborhood of $M$ in
$\bC^{n+1}$.   Let $M'\subset \C^{n'+1}$ be a nondegenerate hyperquadric. 
Suppose that 
$ n' \leq 3( n-e_0(M,p) )$ 
for some point $p\in M$. 
If  $H:U\to \C^{n'+1}$ is a holomorphic mapping with $H(M) \subset
M'$ , then one of the following
mutually exclusive conditions must hold.
\begin{enumerate}
\item [(i)] There is an open subset $V\subset U$ with $M\subset V$ and 
$H(V) \subset
M'$.
\item [(ii)] $H$ is  transversal to $M'$ at $H(p)$ for all $p\in M$
outside a proper real-analytic subset.
\end{enumerate}
\end{Thm}

We should remark that conclusion (ii) in Theorems \ref{main1} and \ref{thm3} cannot
be replaced by the stronger conclusion that  transversality holds for
every $p\in M$, as is shown by Examples \ref{ex1} and \ref{ex2}. In
the equidimensional case, condition \eqref{eigen20} is always
satisfied. The conclusion of Theorem \ref{main1}, in this case, can
be deduced from known results (e.g.\ \cite{brgeom} and \cite{er06})
by using also Theorem
\ref{thm5} of the present paper. Even in the equidimensional case, the conclusion (ii)
cannot be replaced by that of  transversality for all $p\in M$ as is
shown by Example 6.2 in
\cite{er06}. However, if the condition that $M$ is of finite type is added, then it is
unknown if this replacement can be made (see Conjecture 2.7 in \cite{lm06}; see
also Question 1 in \cite{er06}). 

The following result shows that the condition in Theorem \ref{thm3} requiring $M'$
to be a nondegenerate hyperquadric cannot be replaced by the weaker assumption that
$M'$ is a Levi nondegenerate hypersurface. As mentioned above, it also shows that the
condition
\eqref{eigen20} in Theorem \ref{main1} is optimal. 

\begin{Thm}\Label{Thm4} Given $M \subset \C^{n+1}$ a nondegenerate
hyperquadric, there exist a Levi nondegenerate hypersurface $M'\subset \bC^{2n+2}$
and
$H:\C^{n+1}
\to \C^{2n+2}$ a polynomial embedding of degree $2$ such that $H$ sends $M$ into
$M'$, but neither {\rm (i)} nor {\rm (ii)} of Theorem $\ref{main1}$ holds. More
precisely, if $M:=\{Z\in \bC^{n+1}\colon \rho(Z,\bar Z)=0\}$  with 
 \begin{equation}   
\rho(Z,\bar Z):=\im w-\sum_{j=1}^n \delta_j |z_j|^2,\quad Z=(z,w)\in \bC^n\times
\bC,
\end{equation}
 where $\delta_j=\pm 1$ and $\langle \cdot, \cdot\rangle$ is a nondegenerate
hermitian form in $\C^{2n+1}$ with $n$ negative  and $n+1$ positive eigenvalues,
then there exist  a polynomial embedding of degree two, $H:\C^{n+1}
\to \C^{2n+2}$, and a real bihomogeneous polynomial
$\phi(z',\bar {z'})$, $z'\in\C^{2n+1}$, of bidegree $(2,2)$, such that if
\begin{equation}\Label{rho'}\rho'(z',w',\bar z,\bar w):= \im w' -  \langle
z',\bar{z'}\rangle  -\phi(z',\bar{z'})\end{equation} then  
$\rho'\circ H= -4\rho^2$. 
\end{Thm}

We should point out that if the target hypersurface $M'$ in either Theorem
\ref{main1} or \ref{thm3} does not contain any nontrivial complex subvarieties, then
condition (i) is equivalent to the mapping $H$ being constant. Hence, if the
hypothesis that
$M'$ does not contain any nontrivial complex subvarieties is added to either
Theorem
\ref{main1} or \ref{thm3} and $H$ is assumed to be nonconstant, then the
conclusion (ii) necessarily follows. In the last section of this paper, we give a
number of sufficient conditions for (i) to hold (see Theorems \ref{thm5}, \ref{thm6} and 
corollaries). In the equidimensional case (i.e.\
$n=n'$), we give two conditions equivalent to (i) (see Corollary \ref{cor5}).

\section{Examples and a lemma}

In this section, we give some examples, which show that our main results are sharp.
We begin with the following lemma, which expresses conditions (i) and (ii) in
Theorems \ref{main1} and \ref{thm3} in terms of local defining functions for $M$ and
$M'$. 

\begin{Lem}\Label{lem1} Let $M\subset \bC^{n+1}$ and 
$M'\subset \bC^{n'+1}$ be connected
real-analytic hypersurfaces and $U$ an open neighborhood of $M$ in
$\bC^{n+1}$. Let 
 $p\in M$ and $p'\in M'$ and suppose 
that
$M$ and $M'$ are defined locally by $\rho=0$ and $\rho'=0$ near $p$ and $p'$,
respectively. Let $H\colon U\to \bC^{n'+1}$ be
a holomorphic mapping with $H(M)\subset M'$ and $H(p)=p'$. If {\rm
(i)} in Theorem $\ref{main1}$ does not hold, then there exists a unique integer $k\geq
1$ such that $\rho'\circ H=a\rho^k$, where $a$ is a real-valued, real-analytic function
defined near $p$ in $\bC^{n+1}$ with $a|_M\not\equiv 0$. Moreover, the condition {\rm
(ii)} in Theorem $\ref{main1}$ is equivalent to $k=1$. 
\end{Lem}

\begin{proof}  Since $H(M)\subset M'$ and (i) does not hold, $\rho'\circ H$
vanishes on $M$ but is not identically zero near $p$ (by Remark \ref{rmk1}).
Hence,
$\rho'\circ H=b\rho$, where $b\not\equiv 0$. By unique factorization, there is a
unique integer $l\geq0$ such that $b=a\rho^l$ with $a|_M\not\equiv 0$. Now, (ii)
is equivalent to
$b|_M\not \equiv 0$, in view of Remark \ref{rmk1}, and hence $k=1+l=1$. The
conclusion of the lemma now follows. 
\end{proof}

\begin{Ex}\Label{ex0}{\rm Let $M \subset \C^2$ be the Levi-flat hypersurface given
by
$\rho(z,w,\1z,\1w) :=\im w=0$, and $M' \subset \C^2$ the hypersurface given by 
$\rho'(z',w',\1z',\1w') :=\im w' - |z'|^2=0$.  The mapping $H: \C^2 \to\C^2$ given by
$H(z,w)=(w,iw^2)$ maps $M$ into $M'$, but satisfies neither (i) nor (ii) of
Theorem \ref{main1}.  Indeed, we have $\rho'\circ H \equiv -2\rho^2$. 
Note that $M$ is not
holomorphically nondegenerate, which is the only assumption of Theorem
\ref{main1} in this case ($n'=n=1$).}
\end{Ex}

\begin{Ex}\Label{ex1} {\rm Let $M\subset \bC^3$ be the unit sphere,
$$
\rho(Z,\bar Z):=|Z_1|^2+|Z_2|^2+|Z_3|^2-1=0,
$$
and $M'\subset \bC^5$ be the hyperquadric defined by 
$$
\rho'(z,',w',\bar z',\bar w'):=\im w'-(|z_1'|^2+|z_2'|^2+|z_3'|^2-|z_4'|^2)=0.
$$
Consider the mapping
$$
H(Z):=(Z_1Z_2,Z_2^2,Z_2Z_3,Z_2,0).
$$
Observe that we have the identity
$\rho'(H(Z),\overline{H(Z)})=-|Z_2|^2\rho(Z,\bar Z)$. We conclude that $H$
sends
$M$ into $M'$, $H(\bC^3)$ is not contained in
$M'$, and
$H$ is not  transversal to $M'$ at $0$. Note also that Theorem
\ref{main1} applies, since condition \eqref{eigen10} holds. This
example shows that, under assumption \eqref{eigen10},  conclusion
(ii) in Theorem
\ref{main1}  cannot be replaced by the stronger
conclusion of  transversality at all points of
$M$. 
 }
\end{Ex}

\begin{Ex}\Label{ex2}{\rm Let $M \subset \C^2$ be the hypersurface given by
$\rho(z,w,\1z,\1w) :=\im w - |z|^2=0$, $M'\subset \C^3$ the Levi-nondegenerate
hyperquadric given by $\rho'(z',w',\1z',\1w'): = \im w'+|z'_1|^2-|z'_2|^2 = 0$ and
$H:\C^2\to\C^3$ given by
$$H(z,w) = (z+z^2 +\frac{i}{2}w, z-z^2 -\frac{i}{2}w, -2zw).
$$ We have $\rho'\circ H \equiv -2(z+\1z)\rho$.  Hence $H$ is 
transversal on
$M$ outside the real-analytic submanifold of $M$ given by $\re z=0$.  For every $p'
\in M'$, we have $e_0(M',p') = 0$ and $e(M',p') = 1$.  Hence,
\eqref{eigen20} in Theorem \ref{main1} holds (but
\eqref{eigen10} does not).  Also, the assumption on $n'$ in Theorem \ref{thm3} holds. 
Moreover, $M$ is holomorphically nondegenerate, (i) does
not hold, and   transversality does not hold at every point of $M$.
This example shows that (ii) in Theorems \ref{main1} and
\ref{thm3} cannot be replaced by the stronger condition of 
transversality at all points of
$M$. Observe that $H(\bC^2)$ is the 2-dimensional complex manifold
given by $(z_1'+z'_2)(z'_1-z'_2)+iw'=(z'_1+z'_2)^3/2$. 
}
\end{Ex}

\begin{Ex}\Label{ex3}{\rm Let $M \subset \C^2$ be the hypersurface given by
$\rho(z,w,\1z,\1w) :=\im w - |z|^2=0$, $M'\subset \C^5$ the Levi-nondegenerate
hyperquadric given by $$\rho'(z',w',\1z',\1w'): = \im
w'+|z'_1|^2-|z'_2|^2-|z'_3|^2-|z'_4|^2 = 0,$$ and
$H:\C^2\to\C^5$ given by
$$H(z,w) = (iz+zw,-iz+zw,w,\sqrt 2 z^2, iw^2).
$$ We have $\rho'\circ H \equiv -2 \rho^2$.   Since neither (i) nor (ii)  of Theorem
\ref{thm3} holds, this example shows that the condition $n'\leq 3(n-e_0(M,p))$ cannot
be replaced by $n'\leq 3(n-e_0(M,p))+1$. }
\end{Ex}

\section{Proof of Theorem $\ref{main1}$}

For the proof of Theorem $\ref{main1}$ stated in the introduction, we shall need a
number of preliminary results, which may be of independent interest. Recall that if $M$
is a real-analytic hypersurface in $\bC^{n+1}$, defined locally near $p_0\in M$ by
the real-analytic equation $\rho(Z,\bar Z)=0$, then the {\it Segre variety} 
of $M$ at
$p$, sufficiently close to $p_0$, is given by the holomorphic equation $\rho(Z,\bar
p)=0$. We shall denote the Segre variety of $M$ at $p$ by $\Sigma_p$.  The 
following
proposition will be useful in the proofs of the main results.

\begin{Pro}\Label{pro1} Let $M\subset \C^{n+1}$, and $M'\subset \C^{n'+1}$ be
connected real-analytic hypersurfaces and $U$ an open neighborhood of
$M$ in $ 
\C^{n+1}$.   If $H:U\to \C^{n'+1}$ is a holomorphic mapping with $H(M) \subset
M'$ and $M$ is holomorphically nondegenerate, then at least one of the following
holds.
\begin{enumerate}
\item [(i)] There is an open subset $V\subset U$ with $M\subset V$ and 
$H(V) \subset
M'$.
\item [(iii)] For every $p \in M$ outside a proper real-analytic subset, the rank of
$H|_{\Sigma_p}$  at $p$  is $n$.  
\end{enumerate}
\end{Pro}

\begin{Rem} \Label{rmk2}  {\rm We note that if, for some point
$p\in M$, the restriction of $H$ to the Segre variety of $M$ at $p$ has rank $n$ at
$p$, then (iii) in Proposition \ref{pro1} holds. (This is true even without assuming that
$M$ is holomorphically nondegenerate.) Indeed, the rank at $p$ of $H|_{\Sigma_p} $
is the rank of the
$n\times (n'+1)$ matrix given by $(L_jH_k(p))$, where $j=1,\ldots, n$, $k=1,\ldots,
n'+1$ and $L_1,\ldots, L_n$ is a real-analytic local basis of the $(1,0)$-vector fields
tangent to $M$. Thus, if the rank of this matrix is $n$ at some point, then it is $n$
outside a proper real-analytic subset of $M$ near $p$. A standard connectedness
argument shows that (iii) holds. 
 }
\end{Rem}

\begin{proof}[Proof of Proposition $\ref{pro1}$] We assume, in order to reach a
contradiction, that neither (i) nor (iii) of Proposition \ref{pro1}
holds. Let $p_0\in M$ be a point at which $M$ is finitely
nondegenerate\footnote{Recall
that a real-analytic hypersurface $M\subset \C^{n+1}$ locally
defined near a point $p_0 \in M$ by $\rho(Z,\bar Z) = 0$ is {\it
finitely nondegenerate} at $p_0$ if the  vectors $L^\alpha
\rho_Z(p_0)$, $\alpha \in \mathbb Z_+^n$, span all of $\C^{n+1}$,
where
$\rho_Z$ is the gradient vector of $\rho$ with respect to $Z$ and
$L^\alpha = L_1^{\alpha_1}\ldots L_n^{\alpha_n}$, with
$L_1,\ldots, L_n$  as in Remark \ref{rmk2}.  If
$M$ is connected and holomorphically nondegenerate, then it is
finitely nondegenerate on a dense and open subset of $M$ (see \cite{BHR96} and
\cite{book}).}, and $\rho$, $\rho'$ local defining functions for $M$
and $M'$ near $p_0$ and
$H(p_0)$, respectively. By Lemma \ref{lem1}, there exists an integer
$k\geq 1$ such that 
\begin{equation}\Label{nowCRt0}
\rho'\circ H=a\rho^k,
\end{equation}
where $a$ is not identically zero on $M$. By moving to a nearby point, if necessary,
we may assume that $a(p_0)\neq 0$. We choose normal coordinates $(z,w)\in \bC^n\times\bC$ and
$(z',w')\in \bC^{n'}\times\bC$ for $M$ and $M'$ vanishing at $p_0$ and
$H(p_0)$, respectively. Hence, the defining equations of  $M$ and $M'$ can be
written as $w=Q(z,\bar z,\bar w)$ and
$w'=Q'(z',\bar z',\bar w')$, respectively, with $Q(z,0,\tau)\equiv Q(0,\chi,\tau)\equiv 
Q'(z',0,\tau)\equiv Q'(0,\chi',\tau)\equiv \tau$.
We write $H(z,w)=(f(z,w),g(z,w))$ with $f=(f_1,\ldots,
f_{n'})$. 
It follows from
\eqref{nowCRt0} that
\begin{equation}\Label{basic2}
g(z,w)-Q'(f(z,w),\bar f(\chi,\tau),\bar g(\chi,\tau))=a(z,w,\chi,\tau)(w-Q(z,\chi,\tau))^k,
\end{equation}
with $a(0)\neq 0$. Setting $\chi=0$, $\tau=0$, we have $g(z,w)=a(z,w,0,0)w^k$ and
hence
\begin{equation}\Label{gwk}
g_{w^k}(0)\neq 0.
\end{equation}
We differentiate \eqref{basic2} $k-1$ times with respect to $w$ and then set
$\tau=0$, $w=Q(z,\chi,0)$. We obtain, since $g(z,0)\equiv 0$,
\begin{multline}\Label{gwk-1}
g_{w^{k-1}}(z,Q(z,\chi,0))=\\
\sum_{|\alpha|\leq
k-1}Q'_{(z')^\alpha}(f(z,Q(z,\chi,0),\bar f(\chi,0),0)P_\alpha(f_w(z,Q(z,\chi,0),\ldots,
f_{w^{k-1}}(z,Q(z,\chi,0))),
\end{multline}
where the $P_\alpha(t_1,\ldots,t_{k-1})$ are universal polynomials. We now
differentiate \eqref{gwk-1} with respect to $\chi_j$, for $1\leq j\leq n$, and then
set
$\chi=0$ to obtain
\begin{multline}\Label{chi}
g_{w^k}(z,0) Q_{\chi_j}(z,0,0)=\\
\sum_{|\alpha|\leq
k-1}Q'_{(z')^\alpha\chi'}(f(z,Q(z,\chi,0),\bar
f(\chi,0),0)\bar f_{\chi_j}(0)P_\alpha(f_w(z,Q(z,\chi,0),\ldots,
f_{w^{k-1}}(z,Q(z,\chi,0))).
\end{multline}
Since (iii) does not hold, there exist constants $a_1,\ldots, a_n$ with
$(a_1,\ldots,a_n)\neq 0$ such that
\begin{equation}
\sum_{j=1}^n \bar f_{\chi_j}(0)a_j=0.
\end{equation}
Thus, if we multiply \eqref{chi} by $a_j$ and sum over $j$, then we obtain
\begin{equation}
g_{w^k}(z,0) \sum_{j=1}^na_jQ_{\chi_j}(z,0,0)\equiv 0. 
\end{equation}
It follows from \eqref{gwk} that $\sum_{j=1}^na_jQ_{\chi_j}(z,0,0)\equiv 0$ and,
hence, $\sum_{j=1}^na_jQ_{\chi_jz^\alpha}(0)=0$ for all multi-indices $\alpha$.
This contradicts the finite nondegeneracy of $M$ at $p_0$ and completes the
proof of Proposition \ref{pro1}.  \end{proof}

We mention here that some of the techniques in the proof of Proposition \ref{pro1} were used in \cite{bh}.
The following transversality result may already be known in the folklore. For the
reader's convenience, we include a proof here. 

\begin{Pro}\Label{pro2} Let $M\subset \C^{n+1}$, and $M'\subset \C^{n'+1}$ be
real-analytic hypersurfaces with  
$p\in M$ and $p'\in M'$. Let $H\colon (\bC^{n+1},p)\to (\bC^{n'+1},p')$ be a germ at
$p$ of a holomorphic mapping sending $M$ into $M'$ and such that the restriction
of $H$ to the Segre variety of $M$ at $p$ has rank $n$ at $p$. If
\begin{equation}\Label{eigen1} e(M',p') +e_0(M',p') \le n-1,
\end{equation} 
then either $H$ sends a full neighborhood of $p$ in $\bC^{n+1}$ into $M'$ or
$H$ is  transversal to $M'$ at $p'$. 
\end{Pro}

\begin{proof}[Proof of Proposition $\ref{pro2}$]  We assume, in order to reach a
contradiction, that 
$H$ does not map an open neighborhood of $p$ in $\bC^{n+1}$ into $M'$ and that
$H$ is not  transversal to $M'$ at $p'$.  We choose normal
coordinates $(z,w)\in 
\bC^n\times\bC$ and
$(z',w')\in \bC^{n'}\times\bC$ for $M$ and $M'$, vanishing at $p$
and $p'$, respectively. We write $H(z,w)=(f(z,w),g(z,w))$ with
$f=(f_1,\ldots, f_{n'})$. The defining equations of  $M$ and $M'$
can be written as $w=Q(z,\bar z,\bar w)$ and
$w'=Q'(z',\bar z',\bar w')$, respectively, with $Q(z,0,\tau)\equiv Q(0,\chi,\tau)\equiv 
Q'(z',0,\tau)\equiv Q'(0,\chi',\tau)\equiv \tau$.
The fact that $H$ maps $M$ into $M'$ implies that 
\begin{equation}\Label{basic3}
g(z,w)-Q'(f(z,w),\bar f(\chi,\tau),\bar g(\chi,\tau))=a(z,w,\chi,\tau)(w-Q(z,\chi,\tau)),
\end{equation}
where $a$ is a germ at $0$ of a real-analytic function. Since $H$ is not 
transversal to $M'$ at $p'$, it follows that $a(0)=0$. Let
$v_j:=f_{z_j}(0)$ for $j=1,\ldots, n$. By assumption,
$v_1,\ldots, v_n$ are linearly independent vectors in $\bC^{n'}$. We set
$w=\tau=0$ in
\eqref{basic3}, apply $\partial^2/\partial z_j\partial \chi_l$ and evaluate at $z=\chi=0$
to obtain
\begin{equation}\Label{ortho}
v^*_lA v_j=0, \quad 1\leq j,l\leq n,
\end{equation}
where $A$ is the $n'\times n'$ hermitian matrix $(Q'_{z'_\alpha\chi'_\beta}(0))$, the
vectors $v_j$ are regarded as $n'\times 1$ matrices,  and $^*$ denotes the transpose
conjugate. Note that $A$
represents the Levi form of $M'$ at $p'=H(p)$. We denote the
number of positive, negative, and zero eigenvalues of $A$ by
$e_+,e_-$, and $e_0$, respectively, and observe that $\min
(e_+,e_-)=e(M',p')$ and $e_0=e_0(M',p')$.  Let $E$ be the
$n$-dimensional subspace of $\bC^{n'}$ spanned by
$v_1,\ldots, v_n$. Let $\mathcal L\colon \bC^{n'}\times\bC^{n'}\to
\bC$ be the hermitian form given by $(u,v)\mapsto v^*Au$.  Equation \eqref{ortho} implies that
$\mathcal L$ restricted to $E\times E$ is identically zero. Standard linear algebra
gives $n=\dim E\leq \min(e_+,e_-)+e_0=e(M',p')+e_0(M',p')$,
contradicting
\eqref{eigen1}. This completes the proof of  Proposition \ref{pro2}.
\end{proof}

For the proof of Theorem \ref{main1}, we shall also need the following proposition.

\begin{Pro}\Label{pro3} Let $M\subset \C^{n+1}$, and $M'\subset \C^{n'+1}$ be
connected real-analytic hypersurfaces and $U$ an open neighborhood of
$M$ in $ 
\C^{n+1}$.   
Suppose that 
\begin{equation}\Label{eigen2} n' +e_0(M',p') = 2n, \quad  \
\ \forall p'\in M',
\end{equation}
 holds. Then if  $H:U\to \C^{n'+1}$ is a holomorphic mapping with $H(M) \subset
M'$ such that for every $p\in M$ outside a proper real-analytic subset the restriction
of
$H$ to the Segre variety of $M$ at $p$ has rank $n$ at $p$, then one of the following
mutually exclusive conditions must hold.
\begin{enumerate}
\item [(i)] There is an open subset $V\subset U$ with $M\subset V$ and 
$H(V) \subset
M'$.
\item [(ii)] $H$ is  transversal to $M'$ at $H(p)$ for all $p\in M$
outside a proper real-analytic subset.
\end{enumerate}
\end{Pro}

\begin{proof}[Proof of Proposition $\ref{pro3}$ ] We assume, in order to reach a
contradiction, that neither (i) nor (ii) holds.  Choose $p_0\in M$
such that the restriction of $H$ to $\Sigma_{p_0}$ has rank $n$ at
$p_0$. By Lemma
\ref{lem1}, there exists an integer
$k\geq 2$ such that 
\begin{equation}\Label{nowCRt3}
\rho'\circ H=a\rho^k,
\end{equation}
where $a$ is not identically zero on $M$. By moving to a nearby point, if necessary,
we may assume that $a(p_0)\neq 0$. We choose normal coordinates $(z,w)\in 
\bC^n\times\bC$ and
$(z',w')\in \bC^{n'}\times\bC$ for $M$ and $M'$, vanishing at
$p_0$ and
$H(p_0)$, respectively. We write $H(z,w)=(f(z,w),g(z,w))$ with $f=(f_1,\ldots,
f_{n'})$. As above, the defining equations of  $M$ and $M'$ can be written as
$w=Q(z,\bar z,\bar w)$ and
$w'=Q'(z',\bar z',\bar w')$, respectively.
It follows from
\eqref{nowCRt3} that
\begin{equation}\Label{basic4}
g(z,w)-Q'(f(z,w),\bar f(\chi,\tau),\bar g(\chi,\tau))=a(z,w,\chi,\tau)(w-Q(z,\chi,\tau))^k,
\end{equation}
with $a(0)\neq 0$. Let $v_j:=f_{z_j}(0)$ for $j=1,\ldots, n$. By assumption,
$v_1,\ldots, v_n$ are linearly indepedent vectors in $\bC^{n'}$. 
As in the proof of Proposition \ref{pro2}, we obtain \eqref{ortho}, where $A$ is
as in that proof. We denote the number of zero
eigenvalues of $A$ by $e_0$ and observe that $e_0=e_0(M',p_0')$.  
We introduce the
subspaces $E,F\subset\bC^{n'}$ spanned by $v_1,\ldots, v_n$ and
$Av_1,\ldots, Av_n$, respectively. Observe that the dimension of $F=AE$ is at
least
$n-e_0$. By equation
\eqref{ortho}, it follows that $E$ and $F$ are orthogonal with respect to the
standard hermitian inner product of $\bC^{n'}$ and, hence, $E\cap F=\{0\}$.
Since $n'+e_0 =2n$, we conclude that $\bC^{n'}=E\oplus F$ (and hence the
dimension of $F$ is
$n-e_0$). Let us denote by
$v:=f_w(0)\in \bC^{n'}$. By setting $\chi=0$, $\tau=0$ in \eqref{basic4}, we
conclude that $g(z,w)=a(z,w,0,0)w^k$, and in particular,
$g_w(0)=0$, since $k \ge 2$. By setting
$z=0$,
$\tau=0$ in
\eqref{basic4}, applying 
$\partial^2/\partial
w\partial\chi_j$, for $j=1,\ldots, n$, and evaluating at $0$, we obtain
$v_j^*Av=(Av_j)^*v=0$. Consequently, $v$ is orthogonal to $F$ and, hence,
$v\in E$. We set $z=\chi=0$ in \eqref{basic4}, apply
$\partial^k/\partial w^{k-1}\partial\tau$, and evaluate at $0$. Since $a(0)\neq 0$, we
conclude that
\begin{equation}\Label{contra1}
\bigg(\frac{\partial^{k-1}}{\partial
w^{k-1}}Q'_{\chi'}(f(0,w),0,0)\bigg)\bigg|_{w=0}\bar v\neq 0.
\end{equation}
Similarly, setting $z=0$, $\tau=0$ in \eqref{basic4}, applying $\partial^k/\partial
w^{k-1}\partial\chi_j$, for $j=1,\ldots, n$, and evaluating at $0$, we obtain
\begin{equation}\Label{contra2}
\bigg(\frac{\partial^{k-1}}{\partial
w^{k-1}}Q'_{\chi'}(f(0,w),0,0)\bigg)\bigg|_{w=0}\bar v_j= 0,\quad j=1,\ldots,n.
\end{equation}
Since $v\in E$, \eqref{contra2} contradicts \eqref{contra1},
completing the proof of Proposition  \ref{pro3}.
\end{proof}

We are now in a position to prove Theorem  \ref{main1}. 

\begin{proof}[Proof of Theorem $\ref{main1}$] In view of Proposition \ref{pro1}, we
may assume that (iii) of that proposition holds. If condition \eqref{eigen10} holds,
then conclusion (ii) of Theorem
\ref{main1} follows from Proposition \ref{pro2}. Thus, to complete the proof, we may
assume that condition \eqref{eigen20} holds. Note that if $n' +e_0(M', p_0') <2n$, 
for some $p_0'\in M'$, then $n' +e_0(M', p') <2n$ holds for all $p'$ in an open
neighborhood $V$ of
$p_0'$ in $M'$ since $p'\to e_0(M',p')$ is lower semicontinuous. Moreover, 
condition \eqref{eigen10} then holds for all $p'\in V$.  Indeed, since
$e(M',p')\leq(n'-e_0(M',p'))/2$, it follows that
$e(M',p')+e_0(M',p')\leq (n'+e_0(M',p'))/2<n$ and, hence, \eqref{eigen10} holds in
$V$.  The conclusion of Theorem \ref{main1} follows from Proposition \ref{pro2}
(applied to $V$) and Remark \ref{rmk1}. To complete the proof under condition
\eqref{eigen20}, we may assume that $n'+e_0(M',p')=2n$ for all $p'\in
M'$. The conclusion of the theorem now follows from Proposition \ref{pro3}.
\end{proof}

\section{Proofs of Theorems $\ref{thm3}$ and Theorem $\ref{Thm4}$}

We now give the proofs of Theorems $\ref{thm3}$ and Theorem $\ref{Thm4}$.

\begin{proof}[Proof of Theorem $\ref{thm3}$] We suppose, in order to reach a
contradiction, that (i) and (ii) of Theorem \ref{thm3} both fail. 
Hence, in view of Remark
\ref{rmk1}, we may assume that
$H$ is nowhere  transversal. Let 
$p_0$ be
any point on $M$ at which $e_0(M,p)$ is minimal. Note that $e_0(M,p)$ is then
constant in an open neighborhood of $p_0$. Let us, for brevity, denote
$e_0(M,p_0)$ by $e_0$.  Let $\rho$ and $\rho'$  be
local defining functions for $M$ and $M'$ near $p_0$ and $H(p_0)$, respectively. We
conclude, by Lemma \ref{lem1}, that there exists an integer $k\geq 2$
and a real-analytic function $a$ defined in a neighborhood of
$p_0$, not divisible by $\rho$,  such that 
\begin{equation}\Label{nowCRt}\rho'\circ
H=a\rho^k.
\end{equation}
 Since
$a\not\equiv 0$ on $M$, we may assume, by moving to a nearby point if necessary, 
that $a(p_0)\neq 0$. We choose normal coordinates $(z,w)\in \bC^n\times\bC$ and
$(z',w')\in \bC^{n'}\times\bC$ for $M$ and $M'$, vanishing at
$p_0$ and
$H(p_0)$, respectively, and  write $H(z,w)=(f(z,w),g(z,w))$ with $f=(f_1,\ldots,
f_{n'})$. The defining equation of $M'$ can be written as
\begin{equation}
2i\rho'(z',w',\bar z',\bar w')=w-\bar w-2i\langle z',\bar z'\rangle,
\end{equation}
where $\langle\cdot,\cdot\rangle$ is a nondegenerate hermitian form,
and that of $M$ as
\begin{equation}
2i\rho(z,w,\bar z,\bar w)=w-\bar w-2i\sum_{j=1}^{n-e_0}\delta_jz_j\bar
z_j+\psi(z,\bar z,w+\bar w),
\end{equation}
where $\psi(z,0,s)=\psi(0,\bar z,s)=0$, $\psi(z,\bar z,s)=O (3)$, and
$\delta_j = \pm 1$. Now, identity
\eqref{nowCRt} becomes
\begin{multline}\Label{basic0} g(z,w)-\bar g(\chi,\tau)-2i\langle f(z,w),\bar
f(\chi,\tau)\rangle=\\ b(z,w,\chi,\tau)\bigg(w-\tau
-2i\sum_{j=1}^{n-e_0}\delta_jz_j\chi_j+O(3)\bigg)^k,
\end{multline}
where $b$ is a holomorphic function defined in a neighborhood
of $0$ in $\C^{2n+2}$, with
$b(0)\neq 0$. We introduce the following vectors in $\bC^{n'}$
\begin{equation}\Label{vectors}
v_j:=f_{z_j}(0),\ u_j:=f_{z_j^{k-2}w}(0),\ x_j:=f_{z_jw^{k-1}}(0),\quad j=1,\ldots,
n-e_0,
\end{equation} where the subscripts denote partial derivatives. By
carefully identifying appropriate monomials on both sides in
\eqref{basic0}, we conclude that the following identity of
$3(n-e_0)\times 3(n-e_0)$ matrices holds:
\begin{equation}\Label{bigmatrix0}
\begin{pmatrix}
\langle v_j, \bar v_k\rangle & 
\langle v_j, \bar u_k\rangle & \langle v_j, \bar
x_k\rangle\\ {}&{}&{}\\
\langle u_{j}, \bar v_k\rangle &
\langle u_j,
\bar u_{k}\rangle & \langle u_j, \bar x_k\rangle\\ 
{}&{}&{}\\
\langle x_j, \bar v_k\rangle & 
\langle x_j, \bar u_{k}\rangle & \langle x_j, \bar
x_k\rangle
\end{pmatrix} = 
\begin{pmatrix}
0&0&D_1\\ {}&{}&{}\\ 0&D_2& A_1\\
{}&{}&{}\\
\bar D_1 & \bar A_1 & A_2
\end{pmatrix},
\end{equation}
where $j,k=1,\ldots, n-e_0$ and $D_1$, $D_2$, $A_1$, $A_2$ are
$(n-e_0)\times(n-e_0)$ matrices. Moreover, $D_1$ and $D_2$ are invertible diagonal
matrices.  This proves that the matrix on the left in \eqref{bigmatrix0} is invertible and,
hence, the collection of vectors $v_1,\ldots, v_{n-e_0}, u_1,\ldots, u_{n-e_0},
x_1,\ldots, x_{n-e_0}$, given by \eqref{vectors},  are linearly independent. Since
$n'\leq 3(n-e_0)$ by assumption, we conclude that $n'= 3(n-e_0)$ and the vectors
$v_j,u_j,x_j$, for $j=1,\ldots, n-e_0$, form a basis in $\bC^{n'}$. Let now
$y:=f_{z_1^k}(0)$. Again by careful identification of appropriate
monomials in 
\eqref{basic0}, one can check  that 
$\langle y,y\rangle\neq 0$, but 
$$
\langle y,\bar v_j\rangle=\langle y,\bar u_j\rangle=\langle y,\bar x_j\rangle =0,\quad
j=1,\ldots, n-e_0.
$$
This is clearly a contradiction since the vectors  $v_1,\ldots,
v_{n-e_0}, u_1,\ldots, u_{n-e_0}, x_1,\ldots, x_{n-e_0}$ form  a
basis of $\C^{n'}$. This completes the proof of Theorem \ref{thm3}.
\end{proof}

\begin{proof}[Proof of Theorem $\ref{Thm4}$] We
shall take
$H(z,w)=(f(z,w),g(z,w))$, with $f=(f_1,\ldots, f_{2n+1})$, of the form
\begin{equation}\Label{fandg} f(z,w):=\sum_{j=1}^n z_j v_j +
w v_{n+1}+\sum_{j=1}^nz_jw u_j,\quad g(z,w):=2iw^2,
\end{equation} where $v_1,\ldots, v_{n+1},u_1,\ldots, u_n$ are constant linearly
independent vectors in $\bC^{2n+1}$ to be determined. We claim that the vectors
$v_1,\ldots, v_{n+1},u_1,\ldots, u_n$ and a bihomogeneous polymomial
$\phi(z',\bar z')$ of bidegree
$(2,2)$  can be chosen so that
\begin{multline}\Label{basic1} g(z,w)-\bar g(\chi,\tau)-2i\langle f(z,w),\bar
f(\chi,\tau)\rangle-2i\phi(f(z,w),\bar f(\chi,\tau))=
\\ 2i\bigg(w-\tau
-2i\sum_{j=1}^n\delta_jz_j\chi_j\bigg)^2=\\
2i\bigg(w^2+\tau^2-2w\tau-4i\sum_{j=1}^n\delta_jz_j\chi_jw+
4i\sum_{j=1}^n\delta_jz_j\chi_j\tau
-4\sum_{1\leq j,k\leq
n}\delta_j\delta_kz_jz_k\chi_j\chi_k\bigg).
\end{multline} Let us write $\phi(z',\chi')$ in the form
\begin{equation}
\phi(z',\chi'):=T(z',\chi',z',\chi'),
\end{equation} where $T$ is a multilinear form $\bC^{2n+1}\times
\bC^{2n+1}\times\bC^{2n+1}\times\bC^{2n+1}\to\bC$ with the symmetries
\begin{equation}
\begin{aligned} T(X_1,Y_1,X_2,Y_2)=T(X_2,Y_1,
&X_1,Y_2)=T(X_1,Y_2,X_2,Y_1)\\
\overline{T(X_1,Y_1,X_2,Y_2)} &=T(\bar Y_2,\bar X_2,\bar Y_1,\bar X_1)
\end{aligned}
\end{equation} for any $X_1,X_2,Y_1,Y_2\in \bC^{2n+1}$.

Our goal is to find vectors $v_1,\ldots, v_{n+1},u_1,\ldots, u_n$
forming a basis of $\C^{2n+1}$ and a multilinear form $T$ as
above such that \eqref{basic1} holds.  For this, in view of the
choice of $f$ and $g$ given by \eqref{fandg}, it suffices to
establish the two identities:
\begin{equation}\Label{quadra}
\langle f(z,w),\bar f(\chi,\tau)\rangle =
2w\tau+4i\sum_{j=1}^n\delta_jz_j\chi_jw-
4i\sum_{j=1}^n\delta_jz_j\chi_j\tau,
\end{equation}
\begin{equation}\Label{tensor}
T( f(z,w),\bar f(\chi,\tau),f(z,w),\bar f(\chi,\tau)) = 4\sum_{1\leq j,k\leq
n}\delta_j\delta_kz_jz_k\chi_j\chi_k.
\end{equation}
By carefully identifying
all monomials on both sides in \eqref{quadra}, we
conclude that
\eqref{quadra} holds if and only if the following condition
is satisfied:
\begin{equation}\Label{bigmatrix}
\begin{pmatrix}
\big(\langle v_j, \bar v_k\rangle\big)_{j,k=1}^n & 
\big(\langle v_j, \bar v_{n+1}\rangle\big)_{j=1}^n 
& \big(\langle v_j, \bar
u_k\rangle\big)_{j,k=1}^n\\ {}&{}&{}\\
\big(\langle v_{n+1}, \bar v_k\rangle\big)_{k=1}^n&
\langle v_{n+1},
\bar v_{n+1}\rangle & \big(\langle v_{n+1}, \bar u_k\rangle\big)_{k=1}^n\\ 
{}&{}&{}\\
\big(\langle u_j, \bar v_k\rangle\big)_{j,k=1}^n & 
\big(\langle u_j, \bar v_{n+1}\rangle\big)_{j=1}^n & \big(\langle u_j, \bar
u_k\rangle\big)_{j,k=1}^n
\end{pmatrix} = 
\begin{pmatrix}
0_{n\times n}&0_{n\times 1}&D\\ {}&{}&{}\\ 0_{1\times n} &2& 0_{1\times n}\\
{}&{}&{}\\
\bar D & 0_{n\times 1} & 0_{n\times n}
\end{pmatrix},
\end{equation}
where $D$ is the diagonal $n\times n$ matrix given by 
\begin{equation}
D=\begin{pmatrix} -4i\delta_1& 0& \ldots& 0\\ 0& -4i\delta_2&\ldots&0\\
0& 0&\ddots&0\\ 0& 0& \ldots& -4i\delta_n
\end{pmatrix} .
\end{equation}
We must
show that there is a basis $v_1,\ldots,
v_{n+1},u_1,\ldots, u_n$ of vectors in $\bC^{2n+1}$ such that \eqref{bigmatrix}
holds.  We note
that the eigenvalues of the matrix $\Delta$ on the right in
\eqref{bigmatrix} are $2$ with multiplicity 1, 4 with multiplicity
$n$, and -4 with multiplicity
$n$. Let
$e_1,\ldots,e_{2n+1}$ be the standard basis in  $\bC^{2n+1}$ and $Q$ the matrix of
scalar products $$Q=\left(\langle e_\alpha,\bar
e_\beta\rangle\right)_{\alpha,\beta=1}^{2n+1}.$$
Since $Q$ and $\Delta$ have the same number of positive and negative eigenvalues,
there exists an invertible $(2n+1)\times(2n+1)$ matrix $A$ such that
$\Delta=AQA^*$. If we now let $d_1,\ldots,d_{2n+1}$ be the basis
$$
d_\alpha:=\sum_{\gamma=1}^{2n+1}A_{\alpha\gamma}e_{\gamma},
$$
and set $v_j:=d_j$, $j=1,\ldots,n+1$, and $u_j:=d_{n+1+j}$, $j=1,\ldots, n$, then
$v_1,\ldots, v_{n+1},u_1,\ldots, u_{n}$ is a basis for $\bC^{2n+1}$ that satisfies
\eqref{bigmatrix}.

To determine the multilinear form $T$ to satisfy identity
\eqref{tensor}, we set for $j,k=1,\ldots, n,$
\begin{equation}\Label{tensor1}
T(v_j,\bar v_j,v_k,\bar v_k)=T(v_j,\bar v_k,v_k,\bar
v_j)= \begin{cases} {2\delta_j\delta_k,} & \text{if $j\not= k$,}
\\
{4,} &\text{if $j= k$,}
\end{cases} 
\end{equation}
and
$
T(X,\bar Y,Z,\bar W)=0,
$
for all other choices of $X,Y,Z,W$ among the basis vectors $v_1,\ldots,
v_{n+1},u_1,\ldots, u_n$ of $\bC^{2n+1}$. 
The multilinear mapping $T$ (and hence the bihomogeneous
polynomial
$\phi$) is then uniquely defined by \eqref{tensor1} and the vanishing condition
following that equation. It is then straightforward to check that
\eqref{tensor} is satisfied.  This completes the proof of Theorem
\ref{Thm4}. 

\end{proof}

\section{Sufficient conditions for mapping $\bC^{n+1}$ into the target hypersurface}

In this section, we give a number of sufficient conditions for conclusion (i) in
Theorems \ref{main1} and \ref{thm3} to hold. 

\begin{Thm}\Label{thm5} Let $M\subset \C^{n+1}$, and $M'\subset \C^{n'+1}$ be
connected real-analytic hypersurfaces and $U$ an open neighborhood of $M$
in $\C^{n+1}$.  Assume that $M$ is holomorphically nondegenerate.  If $H\colon
U\to
\bC^{n'+1}$ is a holomorphic mapping such that $H(M)\subset M'$ and the 
rank of $H$ is $\leq n$ at every $p\in M$, then there is an open subset $V\subset U$
with $M\subset V$ and 
$H(V) \subset
M'$. 
\end{Thm}

\begin{proof} 
We observe that the real rank of $H|_M$ is $\leq 2n$, by assumption. 
We consider first the case where the rank of $H|_M$ is equal to $2n$ at some
point, and hence on an open and dense subset of $M$. Let $p\in M$ be such a point.
We identify $\bC^{n+1}$ with $\bR^{2n+2}$ and denote by $h\colon U\to
\bR^{2n'+2}$ the corresponding real-analytic mapping induced by $H$. Note that the
rank of $h$ at $p$ is $2n$, as is the rank of $h|_M$. It follows that $\ker dh(p)$,
which is a 2-dimensional subspace of $T_p\bR^{2n+2}$, is not contained in the
hyperplane $T_pM\subset T_p\bR^{2n+2}$ and, hence, $\ker dh(p)$ is transversal
to
$T_pM$. Consequently, we can find a $2n$-dimensional submanifold $S\subset M$
through $p$ that is transversal to $\ker dh(p)$ at $p$. By the rank theorem, there is
an open neighborhood $W\subset \bR^{2n+2}$ of $p$ such that
$h(W)=h(S\cap W)$ and, hence, in particular $h(W)\subset h(M\cap
W)\subset M'$. The conclusion of the theorem now follows from
Remark \ref{rmk1}.

To complete the proof, we consider now the case where the rank of $H|_M$ is $\leq
2n-1$ at every point of $M$.  Choose $p_0\in M$ such that $M$ is
finitely nondegenerate at $p_0$ and $H|_M$ has maximal rank $m
\leq 2n-1$ at
$p_0$. (This is possible since the points at which $M$ is finitely
nondegenerate are dense in $M$.)  Let
$\omega$ be a small neighborhood of
$p_0$ in
$M$ such that the rank of
$H|_M$ is constant in $\omega$. The image $H(\omega)$ is then a real-analytic
submanifold of $\bC^{n'+1}$. By moving the point $p_0$ slightly and shrinking
$\omega$ if necessary, we may assume that $H(\omega)$ is a CR submanifold.
Since the rank of $H|_M$ in $\omega$ is $m \leq 2n-1$,
$H(\omega)$ is CR diffeomeorphic to an $m$-dimensional CR submanifold of
$\omega\subset M$. Since the CR dimension $k$ of $H(\omega)$ is $\leq m/2$, it
follows that
$k<n$. In particular, the restriction of $H$ to the Segre variety of $\omega$ at any point
in
$\omega$ has rank $k<n$. The conclusion of the theorem now follows from
Proposition \ref{pro1}, as well as Remark \ref{rmk1}.
\end{proof}

We remark that 
if $M\subset
\bC^{n+1}$ and $M'\subset \bC^{n'+1}$ are connected real-analytic hypersurfaces,
$U$ an open neighborhood of $M$ in $\bC^{n+1}$ and $H\colon U\to
\bC^{n'+1}$ a holomorphic mapping with $H(M)\subset M'$, then (i) in
Theorems \ref{main1} and \ref{thm3} holds if
and only if for every $p\in M$ there is an open neighborhod $W$ of $p$ in
$\bC^{n+1}$ such $H(W)\subset
\Sigma'_{H(p)}$ . Indeed, if $\rho'(Z',\bar Z')=0$ is a defining equation for $M'$ near
$H(p)$, then $H$ sends a full neighborhood of $p$ in $\bC^{n+1}$ into $M'$ if and
only if $\rho'(H(Z),\overline{H(Z)})\equiv 0$. On the other hand, $H$ sends a full
neighborhood of $p$ into $\Sigma'_{H(p)}$ if and only if
$\rho'(H(Z),\overline{H(p)})\equiv 0$. The conclusion above follows easily from these
facts. In general, it is not enough to have an open neighborhod $W$ of a single 
point  $p\in M$ such $H(W)\subset
\Sigma'_{H(p)}$ to conclude that $H$ sends a full neighborhood of $M$ in
$\bC^{n+1}$ into $M'$ as Example \ref{ex2.4} below illustrates. However, it does
suffice in the equidimensional case, as the following straightforward corollary of
Theorem
\ref{thm5} shows. 

\begin{Cor}\Label{cor5} Let $M, M'\subset \C^{n+1}$ be
connected real-analytic hypersurfaces and $U$ an open neighborhood of $M$
in $\C^{n+1}$.  Assume that $M$ is holomorphically nondegenerate and let
 $H\colon U\to
\bC^{n+1}$be a holomorphic mapping such that $H(M)\subset M'$. The following are
equivalent:
\begin{enumerate}
\item [(i)] There is an open subset $V\subset U$ with $M\subset V$ and 
$H(V) \subset
M'$.
\item [(iv)] There exist $p\in M$ and an open neighborhood
$W\subset U$ of $p$ in $\C^{n+1}$ such that $H(M\cap W)$ is
contained in the Segre variety of $M'$ at $H(p)$.
\item [(v)] The rank of $H$ is $\leq n$ at every $p\in M$.
\end{enumerate}
\end{Cor}

Another straightforward corollary of Theorem
\ref{thm5} is the following. 

\begin{Cor}\Label{cor4} Let $M\subset \bC^{n+1}$ be a real-analytic
holomorphically nondegenerate hypersurface with $p\in M$, $M'\subset \bC^N$ a
real-analytic submanifold with $p'\in M'$, and $H\colon (\bC^{n+1},p)\to (\bC^N,p')$ a
germ at 
$p$ of a holomorphic mapping. If there exists a complex subvariety $X\subset \bC^N$
of dimension $\leq n$ and with $p'\in X$, such that $H(M)\subset X\cap M'$, then
$H(\bC^{n+1})\subset M'$.\end{Cor}

\begin{Ex}\Label{ex2.4}{\rm Let $M\subset \bC^2$ be the unit sphere,
$$
\rho(Z,\bar Z):=|Z_1|^2+|Z_2|^2-1=0,
$$
and $M'\subset \bC^3$ the hypersurface defined by 
$$
\rho'(z,',w',\bar z',\bar w'):=\im w'-|z_1'|^2(|z_1'|^2+|z_2'|^2-1)=0.
$$
Consider the mapping
$$
H(Z):=(Z_1,Z_2,0).
$$
Observe that we have the identity $\rho'(H(Z),\overline{H(Z)})=-|Z_1|^2\rho(Z,\bar
Z)$. We conclude that $H$ sends $M$ into $M'\cap \Sigma'_0$, where $\Sigma'_0$
denotes the Segre variety of $M'$ at $0$ and hence (iv) holds.
However, neither (i) nor (v) of Corollary \ref{cor5} holds.
This example shows that the
conclusion of Corollary \ref{cor5} fails if
$M'$ is a hypersurface in $\bC^{n+2}$ rather than $\bC^{n+1}$ as in that
corollary. }\end{Ex}

However, the following result shows that if $M'$ is a nondegenerate hyperquadric
in
$\bC^{n+2}$, then the conclusion of Corollary \ref{cor5} still holds. 

\begin{Cor}\Label{cor6} Let $M\subset \C^{n+1}$ be a
connected, real-analytic, holomorphically nondegenerate hypersurfaces and $U$
an open neighborhood of
$M$ in $\C^{n+1}$.  Let $M'\subset \bC^{n+2}$ be a nondegenerate hyperquadric.  
If $H\colon U\to
\bC^{n+2}$ is a holomorphic mapping such that $H(M)$ is contained in
the intersection of $M'$ with one of its Segre varieties, then there is an open subset $V\subset U$
with $M\subset V$ and 
$H(V) \subset
M'$.
\end{Cor}

\begin{proof} Let $M'$ be given by 
\begin{equation}   
\im w'=\sum_{j=1}^{n+1}\delta_j |z'_j|^2,\quad (z,w)\in \bC^{n+1}\times
\bC,
\end{equation}
where $\delta_j=\pm 1$. Without loss of generality, we may assume that $H(M)$ is
contained in the intersection of $M'$ with the Segre variety of $M'$ at $0$, i.e.\ in the
variety  given by $w'=0$ and $\sum_{j=1}^{n+1} \delta_j |z'_j|^2=0$. Hence,
$H=(f_1,\ldots, f_{n+1},0)$ where
\begin{equation}
\sum_{j=1}^{n+1} \delta_j |f_j(Z)|^2=0,\quad Z\in M.
\end{equation}
We may assume that $H$ is not constant and hence, after
reordering the coordinates if necessary, that $f_{n+1}$ is not
identically 0 on
$M$. For
$Z\in M$ outside the zero set of $f_{n+1}$, we then have 
\begin{equation}\Label{id}
\sum_{j=1}^{n} \delta_j |\tilde f_j(Z)|^2=-\delta_{n+1},
\end{equation}
where $\tilde f_j:=f_j/f_{n+1}$. Hence, the mapping $\tilde H:=(\tilde f_1,\ldots,\tilde
f_n)$ sends $M$, outside the zeros of $f_{n+1}$, into the hyperquadric 
given by $\sum_{j=1}^{n} \delta_j |z'_j|^2=-\delta_{n+1}$ in $\bC^n$. By Corollary
\ref{cor4}, it follows that \eqref{id} holds identically and, hence, 
$H$ sends a neighborhood of $M$ in $\C^{n+1}$ 
into $M'$. 
\end{proof}

\begin{Ex} {\rm Let $M\subset \bC^2$ be the unit sphere,
$$
\rho(Z,\bar Z):=|Z_1|^2+|Z_2|^2-1=0,
$$
and $M'\subset \bC^4$ the hyperquadric defined by 
$$
\rho'(z,',w',\bar z',\bar w'):=\im w'-(|z_1'|^2+|z_2'|^2-|z_3'|^2)=0.
$$
Consider the mapping
$$
H(z):=(Z_1Z_2,Z_2^2,Z_2,0).
$$
Observe that we have the identity
$\rho'(H(Z),\overline{H(Z)})=-|Z_2|^2\rho(Z,\bar Z)$. We conclude that $H$
sends $M$ into $M'\cap \Sigma'_0$, where $\Sigma'_0$ denotes the Segre variety
of $M'$ at $0$. Observe that the dimension of $\Sigma'_0$ is 3 and the CR
dimension of $M$ is 1. Moreover, $H(\bC^2)$ is not contained in $M'$. This
example shows that the conclusion of Corollary  \ref{cor6} fails if $M'$ is a
hyperquadric in $\bC^{n+3}$ instead of $\bC^{n+2}$.
 }
\end{Ex}

We conclude this paper by giving another sufficient condition for (i) in Theorem
\ref{main1} to hold. We should point out that a proof of Theorem \ref{thm6} below  in the case where $M$ and $M'$ are nondegenerate hyperquadrics was given in \cite{bh}.  We use the notation  $e(M,p)$ and $e_0(M,p)$ introduced
in the introduction.

\begin{Thm}\Label{thm6} Let $M\subset \C^{n+1}$, and $M'\subset
\C^{n'+1}$ be connected real-analytic hypersurfaces and $U$ an open
neighborhood of $M$ in $\C^{n+1}$.  Assume that $M$ is holomorphically
nondegenerate and 
\begin{equation}\Label{thm6cond}
e(M',p')+e_0(M',p')<\sup_{q\in M} e(M,q),\quad \forall p'\in M'.
\end{equation}  
If $H\colon U\to
\bC^{n'+1}$ is a holomorphic mapping such that $H(M)\subset M'$, then there
is an open subset $V\subset U$ with $M\subset V$ and 
$H(V) \subset
M'$. 
\end{Thm}

\begin{proof} We first observe that \eqref{eigen10} follows from
\eqref{thm6cond}, since $e(M,q)\leq n/2$ for all $q\in M$. It follows from
Theorem \ref{main1} that either (i) or (ii) of that theorem must hold. Thus, to
complete the proof of Theorem \ref{thm6}, it suffices to show that (ii) cannot
hold. Let us assume, in order to reach a contradiction, that (ii) holds. We note
that $e(M,p)$  is an integer-valued lower semicontinuous function on
$M$. It follows that $e(M,p)=\sup_{q\in M} e(M,q)$ for
$p$ in an non-empty open subset of
$M$. Hence, we can find a point $p\in M$ such that $H$ is tranversal to $M'$
at $p':=H(p)$ and 
\begin{equation}\Label{last}
e(M',p')+e_0(M',p')<e(M,p)=\sup_{q\in M}
e(M,q). 
\end{equation}
We choose normal
coordinates $(z,w)\in 
\bC^n\times\bC$ and
$(z',w')\in \bC^{n'}\times\bC$ for $M$ and $M'$, vanishing at $p$
and $p'$, respectively. We write $H(z,w)=(f(z,w),g(z,w))$ with
$f=(f_1,\ldots, f_{n'})$. The defining equations of  $M$ and $M'$
can be written as $w=Q(z,\bar z,\bar w)$ and
$w'=Q'(z',\bar z',\bar w')$, respectively, with $Q(z,0,\tau)\equiv Q(0,\chi,\tau)\equiv 
Q'(z',0,\tau)\equiv Q'(0,\chi',\tau)\equiv \tau$.
The fact that $H$ maps $M$ into $M'$ implies that 
\eqref{basic3}
in the proof of Proposition \ref{pro2} holds. Since $H$ is 
transversal to $M'$ at $p'$, it follows that $a(0)\neq 0$ (see Remark \ref{rmk1}).
Let $v_j:=f_{z_j}(0)$ for $j=1,\ldots, n$. We let $B$ denote the $n'\times n$
matrix whose columns are $v_1,\ldots, v_n$. We set
$w=\tau=0$ in
\eqref{basic3}, apply $\partial^2/\partial z_j\partial \chi_l$ to both sides of
\eqref{basic3}, for
$1\leq j,l\leq n$, and evaluate at
$z=\chi=0$ to obtain
\begin{equation}\Label{leviid}
B^*A' B=a(0)A,
\end{equation}
where $A'$ is the $n'\times n'$ hermitian matrix
$(Q'_{z'_i\chi'_j}(0))$,  $A$ is the $n\times n$ hermitian matrix
$(Q_{z_k\chi_l}(0))$, and $^*$ denotes the transpose conjugate. Note
that $A$ and $A'$ represent the Levi forms of $M$ and $M'$ at $p$
and $p'$, respectively. Denote by $e_+,e_-,e_0$ and $e'_+,e'_-,e'_0$ the number
of positive, negative, and zero eigenvalues of $A$ and $A'$, respectively. Recall
that $e(M,p)=\min(e_+,e_-)$, $e(M',p')=\min(e'_+,e'_-)$ and
$e_0(M',p')=e'_0$. Thus, the inequality \eqref{last} implies that
$\min(e'_+,e'_-)<\min(e_+,e_-)$, which, by standard linear algebra, 
contradicts identity
\eqref{leviid} with
$a(0)\neq 0$. The proof of Theorem \ref{thm6} is not complete.
\end{proof}

\def\cprime{$'$}


\end{document}